\documentclass[10pt]{amsart}
\usepackage[english]{babel}
\usepackage{amssymb,amsmath,amsthm,amscd
,mathrsfs}
\usepackage{color}   
\newfont{\rams}{msbm10 scaled\magstep1}
\newfont{\ramss}{msbm10 scaled\magstep0}
\newfont{\iams}{msbm10}
\newfont{\gotic}{eufm10 scaled\magstep1}
\newfont{\bellap}{eusm10 scaled\magstep1}

\newcommand{\sss}{{\rm Spec}}

\newcommand{\hgt}{{\rm ht}}
\newcommand{\mmm}{{\rm Max}}

\newcommand{\zar}{{\rm Zar}}
\newcommand{\kr}{{\rm Kr}}
\newcommand{\ff}{\mathscr F}

\newcommand{\ms}{\mathscr}
\newcommand{\stf}{\star{_{_{\!}f}}}
\newcommand{\ad}{\mbox{\it\texttt{Cl}}}
\newcommand{\FF}{\overline{\boldsymbol{F}}}
\newcommand{\cc}{{\boldsymbol{c}}}
\newcommand{\f}{{\boldsymbol{f}}}

\newcommand{\z}{{\ldots}}
\newcommand{\w}{{\setminus}}

\swapnumbers
\newtheoremstyle{break}
  {9pt}
  {9pt}
  {\itshape}
  {}
  {\textsc}
  {.}
  {.7em}
  {}
\newtheoremstyle{break1}
  {9pt}
  {9pt}
  {}
  {}
  {\textsc}
  {.}
  {.7em}
  {}
\theoremstyle{break}
\newtheorem{thm}{ \textsc{Theorem}}[section]

\newtheorem{cor}[thm]{ \textsc{Corollary}}
\newtheorem{lem}[thm]{ \textsc{Lemma}}
\newtheorem{prop}[thm]{ \textsc{Proposition}}
\theoremstyle{break1}
\newtheorem{ex}[thm]{ \textsc{Example}}
\newtheorem{oss}[thm]{ \textsc{Remark}}
\theoremstyle{remark}

\title[Constructible topology on spaces of valuation domains]{The constructible topology on spaces of valuation domains}
\author{Carmelo A. Finocchiaro, Marco Fontana, and K. Alan Loper}
\address{C.F. \& M.F., Dipartimento di Matematica, Universit\`a degli Studi
``Roma Tre'', Largo San Leonardo Murialdo, 1,  00146 Rome, Italy.}
\email{carmelo@mat.uniroma3.it }
\email{fontana@mat.uniroma3.it }
\address{K.A.L., Department of Mathematics, Ohio State University, Newark, OH 43055, USA}
\email{lopera@math.ohio-state.edu}

\thanks{\it Acknowledgments. \rm
During the preparation of this paper, the first two authors
were  partially supported by  a research grant PRIN-MiUR}

\date{rev. Aug 10, 2011}

\begin{document}

\maketitle

\hfill{\footnotesize \sl }%

\bigskip

\begin{abstract} We consider pro\-per\-ties and applications of a compact, Hausdorff topology
 called the ``ultrafilter topology''  defined on an {\sl arbitrary spectral space} and we observe that this topology coincides with the constructible topology. If $K$ is a field and  $A$  a subring of $K$, we show that the space $\zar(K|A)$ of all valuation domains,  having $K$ as quotient field and containing $A$, (endowed with the Zariski topology) is a spectral space by giving in this general setting the explicit construction of a ring whose Zariski spectrum is homeomorphic to $\zar(K|A)$. We  extend results  regarding spectral topologies on the spaces of all valuation domains and apply the theory developed to study  representations of integrally closed domains as intersections of va\-lua\-tion overrings.  As a very particular case, we prove that two collections of valuation domains of $K$ with the same ultrafilter closure represent, as an interse\-ction, the same integrally closed domain.
\end{abstract}

\medskip


\section{Introduction}
{The motivations for studying spaces of valuation domains come from various directions and, historically, mainly from Zariski's work
for building up algebraic geo\-metry by algebraic means (see\cite{za:39} and \cite{zasa}), from rigid algebraic geometry started by J. Tate (see \cite{ta}, \cite{f-vdP}, and \cite{hu-2}) and from real algebraic geometry (see \cite{sch} and \cite{hu-2}); for a deeper insight  on this topics see the paper by Huber-Knebusch \cite{hu-kn}.

Let $K$ be a field and let $A$ be a subring of $K$.  The goal of this paper is to extend results in the literature concerning topologies on the collection of valuation domains which have $K$ as quotient field, and which have $A$ as a subring { and to provide some applications of these results to the representations of integrally closed domains as intersections of valuation overrings}.  We denote this collection by $\zar(K|A)$.  In case $A$ is the prime subring of $K$, then $\zar(K|A)$ includes all valuation domains with $K$ as quotient field and we denote it by simply 
Zar$(K)$.  A first topological approach to the space Zar$(K)$ is due to Zariski that proved the quasi-compactness of this space, endowed with what is now called the Zariski topology (see \cite{za} and \cite{zasa}).
Later, it was proven, and rediscovered by several authors with a variety of different techniques, that if $K$ is the quotient field of $A$ then $\zar(K|A)$ endowed with Zariski's topology  is a spectral space in the sense of Hochster \cite{ho} (see \cite{dofefo}, \cite{dofo}, \cite{hu-kn} and the appendix of \cite{ku}).

  In Section 2, we start by recalling the definition and the basic properties of the constructible topology on an arbitrary topological space, using the notation introduced in  \cite[Section 2]{sch-tr}  (for further information cf.  \cite[\S 4]{ch}, {\cite[({\bf I}.7.2.11) and ({\bf I}.7.2.12)]{EGA}}, \cite {ho}). Then,  we provide a description of the closure in the contructible topology of any subset of a spectral space by using ultrafilters and ``ultrafilter limit points'' (definition given later). As an application, we obtain a new proof that  the ultrafilter topology on the prime spectrum of  commutative ring $R$, introduced in \cite{folo}, is identical to the classical constructible topology  on this space.
  
     Section 3 is devoted to the study of  the space $\zar(K|A)$ for {\sl any} subring $A$ of $K$, endowed  with the Zariski topology or the constructible topology.
    The versatility of the ultrafilter approach to the constructible topology  is demonstrated in this section, and in the following Section 4, where we make use of  Kronecker function rings. 
     The key result in Section 3 is a proof that the space $\zar(K|A)$ is spectral with respect to the constructible {(and to the Zariski)} topology by giving, in this general setting, the explicit construction of a ring whose prime spectrum is canonically homeomorphic to $\zar(K|A)$.
        This is broader than  the results of Dobbs, Fedder, and Fontana (cf.  \cite{dofefo} and \cite{dofo}), who proved their results   in the case  where $K$ is the quotient field of $A$ { (and only considering the case of the Zariski topology)}. 
        
       Especially noteworthy in Section 4  are the applications of the topological properties of $\zar(K|A)$, endowed with the constructible topology (or, with the inverse topology, in the sense of Hochster \cite{ho}), to the representations of  integrally closed domains as intersections of valuation overrings.  
  For example, Proposition 4.1 indicates that two collections of valuation domains with the same constructible closure will represent the same domain.  Similarly, Corollary 4.15  indicates how the constructible topological structure of a collection of valuation domains determines the associated finite-type e.a.b. semistar operation.  We also apply these results to the class of vacant domains (those domains which have a unique Kronecker function ring).  In particular, Corollaries 4.10 and 4.11 use the constructible topology to characterize vacant domains.   We then relate closure in the inverse topology to closure in the constructible topology and restate our results concerning e.a.b. semistar operations in terms of the inverse topology.
For some distinguished classes of domains, other important contributions on this circle of ideas were given for instance in \cite{ol-1},\cite{ol-2}, \cite{ol-3}, and \cite{ol-4}.

\medskip

\section{Preliminaries, Spectral Spaces and Ultrafilter Limit Points}
\label{not}
If $X$ is a set, we  denote by $\boldsymbol{\mathscr{B}}(X)$ the collection of all subsets of $X$, and by $\ms B_{\mbox{\Tiny{\texttt{fin}}}}(X)$ the collection of all finite subset of $X$. Moreover, if $\mathscr G$ is a nonempty subset of $\ms B(X)$, then we will { simply} denote by $\bigcap \mathscr G$ (resp. $ \bigcup\ms G $) the set obtained by intersection  (resp. union) of all  the subsets of $X$ belonging to $\mathscr G$, i.e., $\bigcap \mathscr G := \bigcap \{ G \mid G \in  \mathscr G\}$  (resp. $ \bigcup\ms G:=\bigcup\{G\mid G\in \ms G\}$).

 Recall that a nonempty collection $\ff$ of subsets of $X$ is said to be a \emph{filter on}  $X$ if the following conditions are satisfied:
{\rm (a)} $\emptyset \notin \ff$;\ 
{\rm (b)} if $F,G\in \ff$, then $F\cap G\in \ff$;\
{\rm (c)} if $F,G\in \ms B(X)$, $F\subseteq G$, and $F\in \ff$, then $G\in \ff$.

Let $\boldsymbol{\mathcal F}(X)$ be the set of all filters on $X$, partially ordered by inclusion.  We say that a filter $\ff$ on $X$ is an  \emph{ultrafilter on} $X$ if it is a maximal element in $\boldsymbol{\mathcal F}(X)$.  In the following, we will denote the collection of all ultrafilters on a set  $X$ by $\boldsymbol{\beta}(X)$.

For each  $x\in X$,  it is immediately seen that $\beta_X^x:=\beta^x:=\{Z\in\ms B(X)\mid x\in Z\}$ is an ultrafilter on $X$, called the \emph{trivial \emph{(or} fixed \emph{or} principal\emph{)} ultrafilter of $X$ {centered on $x$}}. 

\smallskip

 Recall that a {\it spectral space} is a topological space homeomorphic to the prime spectrum of a ring, equipped with the Zariski topology.  The spectral spaces were characterized by Hochster in 1969 as  quasi-compact Kolmogoroff topological spaces, with a quasi-compact open basis stable under finite intersections  and such that every nonempty irreducible closed subspace has a generic point \cite[Proposition 4]{ho}.

\smallskip

Let $\mathcal X$ be a topological space.  With the notation used in \cite[Section 2]{sch-tr} we set:
$$
\begin{array}{rl}
\mathring{\mathcal K} :=\mathring{\mathcal K}(\mathcal X):=&\hskip -5pt  \{ U \mid U \subseteq \mathcal X,\, U \mbox{ open and quasi-compact in } \mathcal X\},\\
\overline{\mathcal K}:=\overline{\mathcal K}(\mathcal X):=&\hskip -5pt  \{ \mathcal X \setminus U \mid U\in \mathring{\mathcal K}(\mathcal X) \},\\
{\mathcal K} :={\mathcal K}(\mathcal X):=&\hskip -5pt  \mbox{the Boolean algebra of the subsets of $\mathcal X$ generated }\\
&\hskip -5pt \mbox{by $\mathring{\mathcal K}(\mathcal X)$}, \\
\end{array}
$$
 i.e.,  ${\mathcal K}(\mathcal X)$ is the smallest subset of $\mathscr{B}(\mathcal X)$ containing $\mathring{\mathcal K}(\mathcal X)$ and closed with respect to  finite $\cup$\ \!, $\cap$\ \!, and complementation.   As in \cite{sch-tr}, we call the {\it constructible topology on $\mathcal X$} the topology on $\mathcal X$ having ${\mathcal K}(\mathcal X)$ as a basis {(for the open sets). 
We denote by ${\mathcal X}^{\mbox{\tiny{\texttt{cons}}}}$ the set $\mathcal X$ equipped with the constructible topology and we call {\it constructible sets of $\mathcal X$} the elements of ${\mathcal K}(\mathcal X)$  (for Noetherian topological spaces, this notion coincides with that given in \cite[\S 4]{ch}) and {\it proconstructible sets} the closed sets of ${\mathcal X}^{\mbox{\tiny{\texttt{cons}}}}$. 

Now consider $\mathcal Y$  a subset of $\mathcal X$. In the following, we  denote by $\ad(\mathcal Y)$ (respectively, $\ad^{\mbox{\tiny{\texttt{cons}}}}(\mathcal Y)$) the closure of $\mathcal Y$, with respect to the given  topology (respectively, the constructible topology) on $\mathcal X$.

\smallskip

Assume  that $\mathcal X$ is a {\sl spectral} space. In this case, the set 
$\mathring{\mathcal K} \ (:= \mathring{\mathcal K}(\mathcal X)$)
  is a basis of the topology  on $\mathcal X$ and it is closed under finite intersections.  The constructible topology on $\mathcal X$ is the coarsest topology for which $\mathring{\mathcal K}$ is a collection of clopen sets   and $\mathcal X^{\mbox{\tiny{\texttt{cons}}}}$  is a compact, Hausdorff topological space. 
  
  We can consider on  $\mathcal X$ the usual partial order, defined by
$$
x\preccurlyeq y :\Longleftrightarrow y\in\ad(\{x\})\,.
$$
 If $\mathcal Y$ is a subset of $\mathcal X$,  set 
$$
\mathcal Y^{\mbox{\tiny{\texttt{sp}}}}:=\{x\in \mathcal X \mid y\preccurlyeq x, \mbox{for some } y\in \mathcal Y\}\,,\;
\mathcal Y^{\mbox{\tiny{\texttt{gen}}}}:=\{x\in \mathcal  X \mid
x\preccurlyeq y, \mbox{for some } y\in \mathcal Y\}\,.
$$
Then $\mathcal Y^{\mbox{\tiny{\texttt{sp}}}}$ (respectively, $\mathcal Y^{\mbox{\tiny{\texttt{gen}}}}$) is \emph{the closure under specializations} (respectively, \emph{the closure under generizations} or \emph{the generic closure}) of $\mathcal Y$.

Following \cite{ho}, we can  also endow  the spectral space $\mathcal X$ with the so called \emph{inverse topology} (or \emph{dual topology}), that is the topology whose basis of closed sets is the set $\mathring{\mathcal K}(\mathcal X)$ of all open and quasi-compact subspaces of $\mathcal X$ (with respect to the given spectral topology). We  denote by $\mathcal X^{\mbox{\tiny{\texttt{inv}}}}$ the set $\mathcal X$, endowed with the inverse topology.   By \cite[Proposition 8]{ho}, $\mathcal X^{\mbox{\tiny{\texttt{inv}}}}$ is a spectral space and its constructible topology is clearly equal to the the constructible topology associated to the given spectral topology on $\mathcal X$. The following fact provides a motivation for the name given to this topology.

\begin{prop}\label{inv}{\rm(\cite[Proposition 8]{ho})}
 Let $\mathcal X$ be a spectral space.  Denote by $\preccurlyeq$ (respectively, $\preccurlyeq'$) the order induced by the given spectral topology (respectively, the inverse topology)   on $\mathcal X$.  Then, for any $x,\ y\in \mathcal X$, we have 
$$
x\preccurlyeq y \; (\mbox{i.e., }  y\in \ad(\{x\})) \; \Leftrightarrow \;  y\preccurlyeq' x \; (\mbox{i.e., }  x\in \ad^{\mbox{\tiny{\texttt{inv}}}}(\{y\}).
$$
\end{prop}

\begin{oss}\label{inv-spectral} \rm
Let $\mathcal X$ be a spectral space and $\mathcal Y$ be a subset of $\mathcal X$. Then, by  \cite[Lemma 1.1]{fo}, \cite[Corollary to Theorem 1]{ho}  and Proposition \ref{inv} we have
$$
\ad(\mathcal Y)=(\ad^{\mbox{\tiny{\texttt{cons}}}}(\mathcal Y))^{\mbox{\tiny{\texttt{sp}}}} \;\; \mbox{  and  } \;\;
\ad^{\mbox{\tiny{\texttt{inv}}}}(\mathcal Y)=(\ad^{\mbox{\tiny{\texttt{cons}}}}(\mathcal Y))^{\mbox{\tiny{\texttt{gen}}}}\,.
$$
\end{oss}

\smallskip

\begin{prop} 
 \label{ultrafilter-point} 
Let $\mathcal X$ be a spectral space, $\mathcal Y$ be a subset of $\mathcal X$ and $\ms U$ be an ultrafilter on $\mathcal Y$. Set 
$$
\begin{array}{rl}
\overline{\mathcal K}_{\mathcal Y, \ms U} :=& \{\mathcal X\setminus U \mid U\in\mathring{\mathcal K} \mbox{ and } \mathcal Y \setminus U\in \ms U\} \\
=& \{ C \in \overline{\mathcal K} \mid C \cap \mathcal Y \in \ms U \} \\
{\mathcal K}_{\mathcal Y, \ms U} :=& \overline{\mathcal K}_{\mathcal Y, \ms U} \bigcup\  \{U\in \mathring{\mathcal K} \mid  U\cap \mathcal Y\in\ms U\}.
\end{array}
$$
Then, the following statements hold.
\begin{enumerate}
\item[\rm (1)] The set $\bigcap {\mathcal K}_{\mathcal Y, \ms U}$ is a singleton, the set  $\overline{K}_{\mathcal Y, \ms U}:= \bigcap \overline{\mathcal K}_{\mathcal Y, \ms U}$ is an irreducible closed subset of $\mathcal X$. The generic point of $\overline{K}_{\mathcal Y, \ms U}$ is the unique point $x_{\ms U} :=x_{\mathcal Y,\ms U}\in\bigcap 
 \mathcal K_{\mathcal Y, \ms U}$. We will call the point $x_{\ms U} \in \mathcal X$ \emph{the ultrafilter limit point of $\mathcal Y$, with  respect to $\ms U$}.
\item[\rm (2)] $\ad^{\mbox{\tiny{\texttt{cons}}}}(\mathcal Y)=\{x_{\mathcal Y, \ms U} \mid \ms U\in \boldsymbol{\beta}(\mathcal Y) \}$.
\end{enumerate}
\end{prop}
\textsc{Proof}.
(1) By construction ${\mathcal K}_{\mathcal Y, \ms U}$  is a collection of closed subsets of $\mathcal X^{\mbox{\tiny{\texttt{cons}}}}$ with the finite intersection property. 
Thus, $\bigcap {\mathcal K}_{\mathcal Y, \ms U}$ is nonempty. Since $\mathcal X$ is, in particular, a topological space satisfying axiom T$_0$, the conclusion will follow if we show that, if $x\in \bigcap  {\mathcal K}_{\mathcal Y, \ms U}$, 
then $\overline{K}_{\mathcal Y, \ms U}=\ad(\{x\})$. 
Since $ \overline{\mathcal K}_{\mathcal Y, \ms U}\subseteq  {\mathcal K}_{\mathcal Y, \ms U}$ and $\overline{K}_{\mathcal Y, \ms U}$ is closed in $\mathcal X$, it follows that $\ad(\{x\})\subseteq \overline{K}_{\mathcal Y, \ms U}$. 

Conversely, let $z \in \overline{K}_{\mathcal Y, \ms U}$ and let $U$ be an open neighborhood of $z$, with respect to the spectral topology. Without loss of generality, we can assume that $U\in \mathring{\mathcal K}$. We have $\mathcal Y\cap U\in \ms U$ (otherwise, since $\ms U$ is an ultrafilter on $\mathcal Y$, $\mathcal Y\setminus U\in \ms U$, so  $\mathcal X \setminus U\in \overline{\mathcal K}_{\mathcal Y, \ms U}$, and  thus in particular  $z \in \mathcal X\setminus U$: a contradiction), hence $U\in {\mathcal K}_{\mathcal Y, \ms U}$   and $x\in U$, in particular. Then, $z\in \ad(\{x\})$.

(2)  Let $\ms U$ be an ultrafilter on $\mathcal Y$ and let $\Omega$ be an open neighborhood of $x_{\ms U}$, with respect to the constructible topology. Since the collection of all clopen sets of $\mathcal X$ is a basis for the open sets of $\mathcal X^{\mbox{\tiny{\texttt{cons}}}}$, we can assume, without loss of generality, that $\Omega = U\cap (\mathcal X \setminus V)$, for some $U,V\in \mathring{\mathcal K} $. It follows immediately that $U\cap \mathcal Y$ and $\mathcal Y \setminus V$ belong to $\ms U$ 
(otherwise, either $\mathcal X \setminus U$ or $V$ would belong to ${\mathcal K}_{\mathcal Y, \ms U}$, then we would have a contradiction since $ \bigcap {\mathcal K}_{\mathcal Y, \ms U} = \{x_{\ms U}\}$ and $x_{\ms U}\ \in U\cap (\mathcal X \setminus V)$). 
Thus, $\Omega\cap \mathcal Y\in \ms U$ and it is, in particular, nonempty. This proves that $x_{\ms U}\in \ad^{\mbox{\tiny{\texttt{cons}}}}(\mathcal Y)$. 

Conversely, let $x\in \ad^{\mbox{\tiny{\texttt{cons}}}}(\mathcal Y)$. Note that the following collection of sets (subsets of $\mathcal Y$):
$$
\mathcal G:=\{\mathcal Y\cap U\cap (\mathcal X \setminus V) \mid U,V\in  \mathring{\mathcal K},\  x\in U\cap(\mathcal X \setminus V)\}
$$
has the finite intersection property, since $\mathring{\mathcal K}$ is a collection of clopen sets of the compact space $\mathcal X^{\mbox{\tiny{\texttt{cons}}}}$. Pick an ultrafilter $\ms U$ on $\mathcal Y$ such that $\mathcal G\subseteq \ms U$ (Lemma 2.1).  We claim that $x=x_{\ms U}$. To see this, since $\mathcal X$ is a T$_0$ space,  it suffices to show that $x$ and $x_{\ms U}$ have the same set of open neighborhoods in $\mathcal X$  (with respect to the given (spectral) topology). 
Let $U$ be an open and quasi-compact neighborhood of $x$. It follows $\mathcal Y\cap U\in \mathcal G\subseteq \ms U$. Thus $U\in {\mathcal K}_{\mathcal Y, \ms U}$ and, in particular, $x_{\ms U}\in U$.   Conversely, assume, by contradiction, that there is an open and compact neighborhood $U$ of $x_{\ms U}$ such that $x\notin U$. Then, $\mathcal Y\cap (\mathcal X \setminus U)\in \mathcal G\subseteq \ms U$. It follows $\mathcal X \setminus U\in \mathcal {\mathcal K}_{\mathcal Y, \ms U}$ and $x_{\ms U}\in \mathcal X \setminus U$,  a contradiction. \hfill$\Box$

\medskip

We apply the previous result to the prime spectrum of a ring.

\begin{cor}\label{ultrafilter-spec-zar}  Let $R$ be a ring,  $X:= \mbox{\em Spec}(R)$ (equipped with the Zariski topology), $Y$ a subset of $X$ and $\ms U$ an ultrafilter on $Y$. Then, 
 $$
P_{\ms U}  :=P_{Y, \ms U} 
:= \{a \in R  \mid V(a) \cap Y \in \ms U \}
$$
is a  prime ideal which coincides with the ultrafilter limit point $x_{\mathcal Y, \ms U}$ of $X$ defined in Proposition \ref{ultrafilter-point}.

\end{cor}
\textsc{Proof}. By an argument similar to that used in \cite[Lemma 2.4]{calota} (see also \cite{folo}) it can be easily shown that $P_{\ms U}$ is a prime ideal of $R$.
 Let $H\in \mathcal K_{Y,\ms U}$. Then $H\cap Y\in\ms U$ and either $H=\bigcup_{i=1}^n D(f_i) \ (\in \mathring{\mathcal K}(\mbox{Spec}(R)))$, for some $f_1, f_2, \z, f_n\in R$, or $H=V(I) \ (\in \overline{\mathcal K}(\mbox{Spec}(R)))$,  for some finitely generated ideal $I$ of $R$. In the first case, we have $\bigcup_{i=1}^nY\cap D(f_i)\in\ms U$, thus $D(f_j)\cap Y\in\ms U$, for some $j$, $1\leq j \leq n$, and so $P_{Y, \ms U}\in D(f_j)\subseteq H$. In the other case, we have  $V( I)\cap Y\in\ms U$ and, since $I$ is finitely generated. $P_{Y, \ms U}\in V(I)=H$. Recalling that $\bigcap \mathcal K_{Y,\ms U}$ is a singleton, the conclusion follows immediately. \hfill$\Box$

\begin{oss}\label{cons=ultra}
Let $\mathcal X$ be a spectral space and $\mathcal Y$ be a subset of $\mathcal X$. We say that $\mathcal Y$ is \emph{ultrafilter closed} in $\mathcal X$ if $x_{\mathcal Y,\ms U}\in \mathcal Y$, for any ultrafilter $\ms U$ on $\mathcal Y$. By Proposition \ref{ultrafilter-point}, it follows that the collection of all the subsets of $\mathcal X$ that are  ultrafilter closed is the family of  closed sets for a topology $\mathcal X$,  that we call the \emph{ultrafilter topology} of the spectral space $\mathcal X$.  
 If we denote  by $\mathcal X^{\mbox{\tiny{\texttt{ultra}}}}$ the space $\mathcal X$ endowed with the ultrafilter topology, then by Proposition \ref{ultrafilter-point}  we have $\mathcal X^{\mbox{\tiny{\texttt{ultra}}}} =\mathcal X^{\mbox{\tiny{\texttt{cons}}}}$. 
 Therefore, from Proposition \ref{ultrafilter-point}  and Corollary \ref{ultrafilter-spec-zar},    when $\mathcal X$ is the prime spectrum of a commutative ring, we reobtain as a particular case \cite[Theorem 8]{folo}. 
\end{oss}

\begin{prop}\label{genericultra-spectral}
Let $\mathcal X$ be a spectral space and $\mathcal Y$ be a quasi-compact subspace of $\mathcal X$. Then, the generic closure $\mathcal Y^{\mbox{\tiny{\texttt{gen}}}}$ of $\mathcal Y$  in $\mathcal X$ is closed  in $\mathcal X^{\mbox{\tiny{\texttt{cons}}}}$. 
\end{prop}
\noindent
\textsc{Proof}. Preserve the notation of Proposition \ref{ultrafilter-point}, and let $\ms U$ be an ultrafilter on $\mathcal Y^{\mbox{\tiny{\texttt{gen}}}}$.  It is sufficient to show that $x_{\ms U}:=x_{\mathcal Y^{\mbox{\Tiny{\texttt{gen}}}},\ms U}\in \mathcal Y^{\mbox{\tiny{\texttt{gen}}}}$. If not, for each $y\in \mathcal Y$, there is an open and compact  open neighborhood $\Omega_y$ of $y$ such that $x_{\ms U}\notin \Omega_y$. By compactness, the open cover $\{\Omega_y\mid y\in \mathcal Y\}$ of $\mathcal Y$  in $\mathcal X^{\mbox{\tiny{\texttt{cons}}}}$ has a finite subcover, say $\{\Omega_{y_i}\mid i=1,2,\z,n\}$.
 It is easily checked that $\mathcal Y^{\mbox{\tiny{\texttt{gen}}}}\subseteq \bigcup_{i=1}^n\Omega_{y_i}$, 
 i.e.,  $\mathcal Y^{\mbox{\tiny{\texttt{gen}}}}=\bigcup_{i=1}^n (\Omega_{y_i}\cap
 \mathcal Y^{\mbox{\tiny{\texttt{gen}}}}) \in \ms U$. 
  Moreover, since $\ms U$ is an ultrafilter   on $\mathcal Y^{\mbox{\tiny{\texttt{gen}}}}$, $\Omega_{y_i}\cap \mathcal Y^{\mbox{\tiny{\texttt{gen}}}}\in \ms U$, for some $i\in\{1,2, \z,n\}$. Thus, by Proposition \ref{ultrafilter-point}(1), we have $x_{\ms U}\in\Omega_{y_i}$, a contradiction.\hfill$\Box$}

\bigskip

%
\section{The Kronecker function ring (after Halter-Koch) and the Zariski-Riemann surface}
 
Let $K$ be a field and let $A$ be {\sl any} subring of $K$. Denote by $\zar(K|A)$ the set of all the valuation domains having $K$ as quotient field and containing $A$  as a subring.

As  is well known, Zariski \cite{za} (or, \cite[Volume II, Chapter VI, \S 1, page 110]{zasa}) introduced and studied the set $Z :=\zar(K|A)$ together with a topological structure defined by taking, as a basis for the  open sets, the subsets  
$
B_{{\!F}}^{{ Z}} := \{V\in {Z}  \mid V \supseteq F\}
$, for $F$ varying in $ \ms B_{\mbox{\Tiny{\texttt{fin}}}}(K)$, i.e., if $F:=\{x_1, x_2, \dots, x_n\}$, with $x_i \in K$, then 
$$
B_F^{{Z}} = \zar(K|A[x_1, x_2, \dots, x_n]).
$$
This topology is called the \emph{Zariski topology on}  $Z=\zar(K|A)$ and $Z$, equipped with this topology, denoted also later by $  Z^{\mbox{\tiny{\texttt{zar}}}}$, is usually called the \emph{(abstract) Zariski-Riemann surface of $K$ over $A$}.

On the set  $Z=\zar(K|A)$ we can also consider the constructible topology, as defined in the previous section, and we denote, as usual,  $  Z^{\mbox{\tiny{\texttt{cons}}}}$ the space $Z$ endowed with the constructible topology.

 In this section, we  show that both $\zar(K|A)^{\mbox{\tiny{\texttt{zar}}}}$ and $\zar(K|A)^{\mbox{\tiny{\texttt{cons}}}}$ are spectral spaces, by giving in this general setting the explicit construction of a ring whose prime spectrum, equipped with the Zariski topology (respectively, constructible topology), is homeomorphic to $\zar(K|A)^{\mbox{\tiny{\texttt{zar}}}}$ (respectively, $\zar(K|A)^{\mbox{\tiny{\texttt{const}}}}$).
 
 \medskip

Let $K$ be a field and $T$ an indeterminate over $K$. For every $W \in \zar(K(T))$, it is well known that $V:=W \cap K \in  \zar(K)$ \cite[Theorem 19.16(a)]{Gilmer} and conversely, for each $V \in  \zar(K)$, there are infinitely many valuation domains $W$ of $K(T)$ such that $W \cap K = V$, called \emph{extensions of $V$ to $K(T)$} \cite[Proposition 20.11]{Gilmer}. Among the extensions of a valuation $V$ of $K$ to $K(T)$, there is a distinguished one, called \emph{the trivial extension of $V$ to $K(T)$}, which is $V(T) := V[T]_{M[T]}$, where $M$ is the maximal ideal of $V$  \cite[Proposition 18.7]{Gilmer}.

\begin{prop}\label{embedding}
Let $K$ be a field and $T$ an indeterminate over $K$. 
\begin{itemize}
\item[\rm (1)]  The canonical map $\varphi: \zar(K(T))^{\mbox{\tiny{\texttt{zar}}}}\rightarrow  \zar(K)^{\mbox{\tiny{\texttt{zar}}}}$, $W\mapsto W\cap K$, is a continuous surjection. 
\item[\rm (2)] Let $ \zar_0(K(T)) := \{ V(T) \in  \zar(K(T)) \mid V \in  \zar(K)\}$. Then,
$$
\varphi\vert_{\zar_0(K(T))}: \zar_0(K(T))^{\mbox{\tiny{\texttt{zar}}}}\rightarrow  \zar(K)^{\mbox{\tiny{\texttt{zar}}}}
$$ 
is a homeomorphism.
\end{itemize}
\end{prop}
\noindent
\textsc{Proof}. (1) The map $\varphi$ is clearly surjective by the previous remarks. It is {also } a continuous map since, for each finite subset $F$ of $K$ and for each basic open set $B^{\zar(K)}_F$ of  $\zar(K)^{\mbox{\tiny{\texttt{zar}}}}$, it is straightforward to see that $\varphi^{-1}(B^{\zar(K)}_F) =B^{\zar(K(T))}_F$.

(2)  It is obvious that  $\varphi\vert_{\zar_0(K(T))}: \zar_0(K(T))^{\mbox{\tiny{\texttt{zar}}}}\rightarrow  \zar(K)^{\mbox{\tiny{\texttt{zar}}}}$ is a  bijection and, by (1),  is a continuous map. The conclusion will follows if we show that the map $\varphi\vert_{\zar_0(K(T))}$ is also open.
Let $h\in K(T)\w \{0\}$, say $$h:=\dfrac{a_0 +a_1T+\z +a_rT^r}{b_0+b_1 T+\z +b_s T^s},$$
with $a_i$ and $b_j$ in $K$, for $i=0, 1, \z,r$ and $ j=0, 1, \z, s$.
Let  $V(T)$ be  a valutation domain in $ \zar_0(K(T))$, let $v$ be the valuation on $K$ defining $V$  
 and let $ v^\ast$  be the valuation associated to $V(T)$, i.e., $v^\ast(T)=v(1)=0$  and, for each nonzero polynomial $f:=a_0+ a_1T+\z+a_rT^r\in K[T]$,  $v^\ast(f):= \inf\{v(a_i)  \mid i=0, 1, \z,r \}$.

 It is easy to see that $V(T) \in B_h^{\zar(K(T))}$  if and only if { $v^\ast(h) \geq 0$, that is, if and only if }
$$ \inf\{v(a_i) \mid i=0, 1, \z,r\} \geq \inf\{v(b_j) \mid j=0, 1, \z, s\}.$$
Now, {for all $i\in\{0, 1, \z,r\}$ and $j\in\{0, 1, \z,s\}$ such that both $a_i$ and $b_j$ are nonzero, 
set}:
$$
F_{ij}:=\left\{\dfrac{a_i}{b_j},\ \dfrac{a_\lambda}{a_i},\ \dfrac{b_\mu}{b_j} \mid \lambda =0, 1,\z,r,\ \mu=0,1, \z, s\right\}.
$$
Then,  it is not hard to verify 
 that $\varphi(B_h^{\zar_0(K(T))})=\bigcup_{i,j}B_{F_{ij}}^{\zar(K)}$ {  (see also the proof of \cite[Lemma 1]{dofo})}, hence the continuous bijective map
 $\varphi\vert_{\zar_0(K(T))}$ is open, and so a homeomorphism. \hfill$\Box$

\bigskip

Now recall  the following key notion introduced by Halter-Koch \cite[Definition 2.1]{hal}, providing an axiomatic approach to the theory of Kronecker function rings.

Let $K$ be field,  $T$ an indeterminate over $K$, and $R$  a subring of $K(T)$.
 We call $R$  \emph{a $K$--function ring} (after Halter-Koch) if $T$ and $T^{-1}$ belong to $ R$ and,  for each { nonzero} polynomial $f\in K[T]$, $ f(0)\in f(T)R$.

We collect in the next proposition several properties of $K$--function rings that will be useful in the following.

\begin{prop}\label{basic}
Let $K$ be a field, $T$ an indeterminate over $K$ and let  $R$ be a subring of $K(T)$. Assume that $R$  is a $K$--function ring.
\begin{enumerate}
\item[\rm (1)] If $R'$ is a subring of $K(T)$ containing $R$, then $R'$ is also a $K$--function ring.
\item[\rm (2)] If $\ms R$ is a nonempty collection of $K$--function rings (in $K(T)$), then $\bigcap \ms R$ is a $K$--function ring.
\item[\rm (3)] $R$ is a B\'ezout domain with quotient field $K(T)$.
\item[\rm (4)] If $f:=f_0+ f_1T+\z+f_rT^r\in K[T]$, then ${ (f_0, f_1,\z, f_r)R} =fR$.
\item[\rm (5)] For every valuation domain $V$ of $K$, $V(T)$ is a $K$--function ring.
\end{enumerate}
\end{prop}
\noindent
\textsc{Proof}. (1), (2), (3) and (4) were proved in \cite[Remarks at page 47 and Theorem (2.2)]{hal}. To prove (5), observe that, if $v$ is the valuation associated to $V$ and { $ v^\ast$ is the trivial extension of  $v$ to $K(T)$  \cite[page 218]{Gilmer}}, then  $v^\ast(T)=v(1)=0$ { or, equivalently, $T$ is invertible in $V(T)$. Moreover, if $f:=f_0+ f_1T+\z+f_rT^r\in K[T]$,  then $v^\ast(f) \leq v(f_0) = v^\ast(f_0)$,} and so $f(0) = f_0 \in fV(T)$.
\hfill$\Box$

\bigskip

The following fact   provides a slight generalization of \cite[Theorem 2.3]{heu} and its proof is similar to that given by O. Kwegna Heubo, which is  based on the work by Halter-Koch \cite{hal}.

\begin{prop}\label{cruciale}
Let $K$ be a field, $T$ an indeterminate over $K$ and $R$ a $K$--function ring. Then, $\zar(K(T)|R) = \zar_0(K(T)|R)$  (i.e.,  for every valuation domain  $W\in \zar(K(T)|R)$,  $W=(W\cap K)(T)$). 
\end{prop}
\noindent
\textsc{Proof}. Let $W$ be a valutation overring of $R$. First, observe that $V:=W\cap K$ is a valuation ring of $K$ \cite[Theorem 19.16(a)]{Gilmer}. Now, let $v$ be a valuation of $K$ defining $V$ and let  
 $f:=f_0+f_1T+ \z+f_rT^r\in K[T]$, $f \neq 0$. By Proposition \ref{basic}(1),  since $R \subseteq W$, $W$ is a $K$--function ring. Let $w$ be a valuation of $K(T)$ defining $W$.  Since $T$ and $T^{-1}$  belong to $W$, we have $w(T)=0$.  Moreover,  $w\vert_K = v$ and so
 $w(f)\geq \inf\{w(f_i) \mid i=0, 1, \z,r\}=\inf\{v(f_i)\mid i=0, 1, \z,r\}$. 
 
 On the other hand, by Proposition \ref{basic}(4), $fR=f_0R+ f_1+\z+f_rR$, and thus $f_i\in fR$, for every $i=0, 1, \z,r$.
 Since $R\subseteq W$,  we have  $f_i\in fW$ and thus $w(f_i)=v(f_i)\geq w(f)$, for every $i=0, 1, \z,r$.  Therefore,  $w(f)=\inf\{v(f_i) \mid i=0, 1, \z, r\}$. This proves that $w = v^\ast$, and hence  $W$ is the trivial extension of $V$ in $K(T)$.\hfill$\Box$

\begin{prop}\label{char-k-funct}
Let $K$ be a field, $T$ an indeterminate over $K$, and  $R$ a subring of $K(T)$.  Then, the following conditions are equivalent. 
\begin{enumerate}
\item[\rm(i)] $R$ is a $K-$function ring.
\item[\rm(ii)] $R$ is  integrally closed in $K(T)$ and $\zar(K(T)|R)= \zar_0(K(T)|R)$. 
\end{enumerate}
\end{prop} 
\noindent
\textsc{Proof.} (i)$\Rightarrow$(ii) is already known  (Propositions \ref{basic}(3) and \ref{cruciale}). \\
(ii)$\Rightarrow$(i). Since $R$ is integrally closed in  $K(T)$ and $\zar(K(T)|R)=\zar_0(K(T)|R)$, then $R=\bigcap \zar_0(K(T)|R)$. Now, the conclusion is clear, by Proposition \ref{basic}(2, 5). \hfill$\Box$

 
 \medskip
 As a consequence of Propositions \ref{embedding}(2),  \ref{ultrahomeo} and \ref{cruciale}, we deduce immediately the following.
 
\begin{cor}\label{cisiamoquasi1}
Let $K$ be a field, $T$ an indeterminate over $K$ and $R\ (\subseteq K(T))$ a $K$--function ring.  Set $A_R := R \cap K$. Then, the canonical map 
$\varphi: \zar(K(T)|R)\longrightarrow  \zar(K|A_R)$, $W\mapsto W\cap K$, is a  topological embedding, with respect to the Zariski topology.
\end{cor}
\noindent

 As an application of the previous corollary we reobtain in particular   \cite[Corollary 2.2, Proposition 2.7 and Corollary 2.9]{heu}. More precisely,
 
\begin{cor}\label{cisiamoquasi2}
Let $K$ be a field, $A$ {\emph{any}} subring of $K$ and $T$ an indeterminate over $K$. Then,
\begin{itemize}
\item[\rm (1)]  $\kr(K|A):=\bigcap\{V(T)\mid V\in \zar(K|A)\}$  is a $K$--function ring. 
\item[\rm (2)]  The canonical map $\varphi: \zar(K(T)|\kr(K|A))^{\mbox{\tiny{\texttt{zar}}}}\rightarrow  \zar(K|A)^{\mbox{\tiny{\texttt{zar}}}}$, $W \mapsto W\cap K$,  is a homeomorphism.
\item[\rm (3)]  The canonical map $\psi: \sss(\kr(K|A)))^{\mbox{\tiny{\texttt{zar}}}}\rightarrow  \zar(K|A)^{\mbox{\tiny{\texttt{zar}}}}$, $Q \mapsto \kr(K|A)_Q\cap K$,  is a homeomorphism.  In particular, $\zar(K|A)^{\mbox{\tiny{\texttt{zar}}}}$ is a spectral space.
\end{itemize}
\end{cor}

\noindent
\textsc{Proof}. (1) By Proposition \ref{basic}(2 and 5), $\kr(K|A)$ is a $K$--function ring { (in $K(T)$).}

{{(2)}}  Let $R := \kr(K|A)$, then clearly $R \cap K$ coincides with the integral closure $\bar{A}$ of $A$ in $K$, and therefore $\varphi(\zar(K(T)|R))=\zar(K|A)$. Now (2) follows from Corollary \ref{cisiamoquasi1}.

(3) Recall that, if $\mathcal A$ is a Pr\"ufer domain and $\mathcal K$ is the quotient field of $\mathcal A$, by \cite[Proposition 2.2]{dofo},  $\zar(\mathcal K| \mathcal A)^{\mbox{\tiny{\texttt{zar}}}}$ is canonically homeomorphic to $\sss(\mathcal A)^{\mbox{\tiny{\texttt{zar}}}}$ (under the map $\mathcal V  \mapsto M_{\mathcal V} \cap \mathcal A$, where  $M_{\mathcal V}$ is the maximal ideal of the valuation domain $\mathcal V$ of $\mathcal K$ containing $\mathcal A$).
Now, the conclusion follows immediately, since $\kr(K|A)$ is a Pr\"ufer domain with quotient field $K(T)$  (Proposition \ref{basic}(3)). 
 \hfill$\Box$

\begin{oss}  Note that the {noteworthy progress provided by }Corollary \ref{cisiamoquasi2} concerns the case  where $A$ is a proper subfield of $K$. As a matter of fact, if $A$ is an integrally closed domain and $K$ is its quotient field,   statements (2) and (3) of  Corollary \ref{cisiamoquasi2} were already proved in \cite [Theorem 2]{dofo}. 
\end{oss}

 By Corollary \ref{cisiamoquasi2}(3), $\zar(K|A)^{\mbox{\tiny{\texttt{zar}}}}$ is a spectral space.   As we will see next, by  applying Proposition \ref{ultrafilter-point} and Remark \ref{cons=ultra}, we reobtain the main statement of \cite[Theorem 3.4]{fi-fo-lo0}.
  
\begin{cor}\label{ultrafilter-zar} Let $K$ be a field and $A$ a subring of $K$. If $Y$ is a nonempty subset of $\zar(K|A)$ (equipped with the Zariski topology) and $\ms U$ is an ultrafilter on $Y$, then 
 $$A_{\ms U}:= A_{Y, \ms U}:=\{x\in K \mid B_x\cap Y\in \ms U\}$$
 is valuation  domain of $K$ containing $A$   and it coincides with the ultrafilter limit point $x_{Y, \ms U}$ of $\zar(K|A)$ defined as in Proposition \ref{ultrafilter-point}. 
Moreover, $\zar(K|A)^{\mbox{\tiny{\texttt{cons}}}}=\zar(K|A)^{\mbox{\tiny{\texttt{ultra}}}}$.
\end{cor}
\noindent
\textsc{Proof.} Preserve the notation  of Proposition \ref{ultrafilter-point}. By  \cite[Lemma 2.9]{calota} (or \cite[Proposition 3.1]{fi-fo-lo0}), we have $A_{Y,\ms U}\in \zar(K|A)$.   The  conclusion will follow by  Proposition \ref{ultrafilter-point} and  Remark \ref{cons=ultra} if we show that $A_{Y,\ms U}\in \bigcap \mathcal K_{Y,\ms U}$. Suppose $H\in \mathcal K_{Y,\ms U}$. By definition, $H\cap Y\in\ms U$ and either $ H=U $ or $ H=X\setminus U $, for some open and compact  subspace $U$ of $\zar(K|A)^{\mbox{\tiny{\texttt{zar}}}}$. Since $\mathcal  B:=\{B_F\mid F\in\ms B_{\mbox{\Tiny{\texttt{fin}}}}(K)\}$ is a basis of the Zariski topology, $U$ is the union of a finite subfamily $\mathcal B'$ of $\mathcal B$. Thus, if $H=U$ (and so $Y\cap U\in\ms U$), there exists a set $B_F\in \mathcal B'$ such that $B_F\cap Y\in \ms U$. Furthermore, for any element $x\in F$, we have $B_F\cap Y\subseteq B_x\cap Y$, and so $B_x\cap Y\in \ms U$.  By definition, it follows $F\subseteq A_{Y,\ms U}$, that is, $A_{Y,\ms U}\in B_F\subseteq \bigcup \mathcal B'=U$. Now, suppose that $H=X\setminus U$ (thus $Y\setminus U\in \ms U$). We want to show that $A_{Y,\ms U}\in X\setminus U$. Assume, by contradiction, that $A_{Y,\ms U}\in B_F$, for some $B_F\in\mathcal B'$. It follows immediately that $B_F\cap Y\in \ms U$ and, finally, $\emptyset = (B_F\cap Y)\cap(Y\setminus U)\in \ms U$, contradiction. 
\hfill$\Box$


 \begin{prop}\label{ultrahomeo} We preserve the notation  of Proposition \ref{embedding} and, now,  let  $\zar(K(T))$ and $\zar(K)$ be endowed with the  constructible topology. Then, the canonical (surjective) map $\varphi: \zar(K(T))^{\mbox{\tiny{\texttt{cons}}}}\rightarrow  \zar(K)^{\mbox{\tiny{\texttt{cons}}}}$ is continuous and (hence) closed. In particular, $\varphi\vert_{\zar_0(K(T))}$ is a homeo\-morphism of $\zar_0(K(T))^{\mbox{\tiny{\texttt{cons}}}}$ onto $\zar(K)^{\mbox{\tiny{\texttt{cons}}}}$. 
 \end{prop}
 \noindent
 \textsc{Proof.} Let $C$ be a closed subset of $\zar(K)^{\mbox{\tiny{\texttt{cons}}}}$, and let $\ms U$ be an ultrafilter on $\varphi^{-1}(C)$. Set $\varphi':=\varphi|_{\varphi^{-1}(C)}:\varphi^{-1}(C)\longrightarrow C$. Then, by \cite[Lemma 2.1(4)]{fi-fo-lo0}, the set $\ms V:=\ms U_{{\varphi'}}$ is an ultrafilter on $C$. Let $V:= A_{\ms V} \in \zar(K)$ and  $W:= A_{\ms U} \in \zar(K(T))$ { then, by a routine argument, it is easy} to check that $\varphi(W)=V \in C$.   By Proposition \ref{ultrafilter-point}(2), it follows  that $\varphi $ is continuous. Moreover (by \cite[Chapter XI, Theorem 2.1]{Du}), $\varphi$ is closed, since $\zar(K(T))^{\mbox{\tiny{\texttt{cons}}}}$ and $\zar(K)^{\mbox{\tiny{\texttt{cons}}}}$ are compact and Hausdorff \cite[Theorem 3.4(5 and 6)]{fi-fo-lo0}. Finally, the last statement is a consequence of the fact that the restriction of $\varphi$ to $\zar_0(K(T))$ is bijective. \hfill$\Box$

 \begin{oss}\label{infinity} Let $K$ be a field, $T$ an indeterminate over $K$, and set $R_0:=\bigcap\{V(T)\mid V\in \zar(K)\}$. Then, by Propositions \ref{basic}(2) and \ref{cruciale}, it follows that $$\zar_0(K(T))=\zar(K(T)|R_0).$$ In particular, $\zar_0(K(T))$ is a closed subspace of  $\zar(K(T))^{\mbox{\tiny{\texttt{cons}}}}$, in view of \cite[Theorem 3.4(2 and 6)]{fi-fo-lo0}.
 \end{oss}

A ``constructible'' version of  Corollary \ref{cisiamoquasi2}(2 and 3) can also be  easily deduced from  the previous considerations.

\begin{cor}\label{cisiamoquasi3}
Let $K$ be a field, $A$ \emph{any} subring of $K$, $T$ an indeterminate over $K$, and let $\kr(K|A)$ be as in Corollary \ref{cisiamoquasi2}. 
\begin{itemize}
\item[\rm (1)]  The canonical map $\varphi: \zar(K(T)|\kr(K|A))^{\mbox{\tiny{\texttt{cons}}}}\rightarrow  \zar(K|A)^{\mbox{\tiny{\texttt{cons}}}}$, defined by $W \mapsto W\cap K$,  is a homeomorphism.
\item[\rm (2)]  The canonical map $\psi: \sss(\kr(K|A)))^{\mbox{\tiny{\texttt{cons}}}}\rightarrow  \zar(K|A)^{\mbox{\tiny{\texttt{cons}}}}$, defined by $Q \mapsto $ $ \kr(K|A)_Q\cap K$,  is a homeomorphism.  In particular, $\zar(K|A)^{\mbox{\tiny{\texttt{cons}}}}$ is a spectral space  canonically homeomorphic to the prime spectrum of the absolutely flat ring canonically associated to the $K$-function ring  $\kr(K|A)$. 

\end{itemize}
\end{cor}
\noindent
\textsc{Proof}. (1) As observed in the proof of Corollary \ref{cisiamoquasi2}(1), $\kr(K|A)$ is a $K$--function ring, hence this statement follows from    Proposition \ref{ultrahomeo}.

 (2) is a consequence of (1),  \cite[Theorem 3.8 and Remark 3.9]{fi-fo-lo0} and  \cite[Propositions 5 and 6, and Theorem 8]{folo}.
\hfill$\Box$

\section{Some Applications}

 Let $K$ be a field, and let $A$ be  a subring of $K$.  In this section,  we use constantly that on the space $\zar(K|A)$ the contructible topology coincides with the ultrafilter topology  (Remark \ref{cons=ultra}) and we give some applications of the  results of the pre\-vious sections to the representations of integrally closed domains as intersections of valuation overrings.

\begin{prop}\label{closure-intersection}
Let $K$ be a field, $A$ be  a subring of $K$ and $U$ be a subset of $Z :=\zar(K|A)$. 
Let $Y'$ and $Y''$ be two subsets of a given subset $U$ of $Z$ 
and assume that their closures  in $U$, with the subspace topology induced by the  constructible topology of $Z$, coincide, 
i.e., $\ad^{\mbox{\tiny{\texttt{cons}}}}(Y') \cap U =\ad^{\mbox{\tiny{\texttt{cons}}}}(Y'') \cap U$. Then,  $\bigcap\{V'\mid V'\in Y'\}=\bigcap\{V''\mid V''\in Y''\}$. { In particular, for each subset $Y$ of $Z$,  
  $$\bigcap\{V\mid V\in Y\}=\bigcap\{W\mid W\in \ad^{\mbox{\tiny{\texttt{cons}}}}(Y) \}.
  $$}
\end{prop}
\noindent
\textsc{Proof.} Assume, by contradiction, that there is an element  { $x_0\in\bigcap\{V'\mid V'\in Y'\}\setminus\bigcap\{V''\mid V''\in Y''\}$,} and pick a valuation domain $V_0\in Y''$ such that $x_0\notin V_0$. By  \cite[Theorem 3.4(2)]{fi-fo-lo0}, the set $\Omega:=U\setminus B_{x_0}$ is an open subset of $U$, with respect to the subspace topology { induced by $Z^{\mbox{\tiny{\texttt{cons}}}}$,} and it contains $V_0$. Since $V_0\in Y''\subseteq \ad^{\mbox{\tiny{\texttt{cons}}}}(Y'') \cap U = $  $\ad^{\mbox{\tiny{\texttt{cons}}}}(Y') \cap U$ and $V_0 \not\in Y'$, then $\Omega\cap Y'$ is nonempty. This implies that there exists a valuation domain { $V'\in Y'$} such that $x_0\notin V'$, a contradiction. \mbox{    } \hfill $\Box$
\medskip


\begin{oss} \rm Note that the previous proposition is stated  very generally using a ``relative-type'' formulation. 
However, it is clear that, if we take any two subsets $Y'$ and $Y''$ of $Z$, the r\^ole of $U$ can be played by any subset of $Z$ (including $Z$) containing $Y' \cup Y''$.
\end{oss}
\medskip

Let $\Sigma$ be a collection of subrings of a field $K$, having $K$ as quotient field. We say that $\Sigma$ is   {\it locally finite} if each nonzero element of $K$ is noninvertible in at most finitely many of the rings belonging to $\Sigma$.

The following easy result will provide a class of integral domains for which the equality { $\bigcap\{V' \mid V'\in Y'\}=\bigcap\{V''\mid V''\in Y''\}$ does not imply, in general, that $\ad^{\mbox{\tiny{\texttt{cons}}}}(Y')=\ad^{\mbox{\tiny{\texttt{cons}}}}(Y'')$.}

\begin{lem}\label{finitecharclos}
Let $K$ be a field and $A$ be a subring of $K$. If $\Sigma$ is an infinite and locally finite subset of { $Z:=\zar(K|A)$, then $\ad^{\mbox{\tiny{\texttt{cons}}}}(\Sigma)=\Sigma\cup\{K\}$.}
\end{lem}
\noindent
\textsc{Proof.} By Proposition \ref{ultrafilter-point}(2) and \cite[Remark 3.2]{fi-fo-lo0}, it is enough to show that { $K= A_{\ms U} \ (=
\{x\in K \mid B_x\cap \Sigma\in \ms U\})$,} for every nontrivial ultrafilter $\ms U$ on $\Sigma$. 
By {contradiction,} assume that there exists an element $x_0\in K\setminus A_{\ms U}$. Then $\Sigma\setminus B_{x_0}\in \ms U$, and  so it is infinite, since $\ms U$ is nontrivial  {(an ultrafilter containing a finite set is trivial)}. This implies that $x_0$ is { noninvertibile in infinitely many  valuation domains belonging to} $\Sigma$, a contradiction. \hfill$\Box$
\medskip

 As a consequence of the previous lemma, we have that, if an integral domain admits two distinct infinite and locally finite representations as intersection of valu\-tation domains, then the converse of Proposition \ref{closure-intersection} does not hold.  An explicit example is given next.
 
 \begin{ex} Let $k$ be a field and let  $T_1, T_2, T_3$ be three indeterminates.  Let  $B$ be the two-dimensional, local  domain  $k(T_3)[T_1, T_2]_{(T_1, T_2)}$  with maximal ideal  $M_B := (T_1, T_2) k(T_3)[T_1, T_2]_{(T_1, T_2)}$, i.e., $B = k(T_3) + M_B$.  Now, let $V$ be { (the rank 1 discrete)} valuation domain defined by $V := k[T_3]_{(T_3)}$ and let  $A$ be the pullback domain given by 
 
 \begin{center}
 $
 A := V + M_B =   k[T_3]_{(T_3)} + (T_1, T_2) k(T_3)[T_1, T_2]_{(T_1, T_2)}.
 $
 \end{center}
  
 Our goal is to represent $A$ as a locally finite intersection of valuation domains in two different ways.  In fact, we can use one description to generate an infinite number of different such representations.  
 
 First, note that $B$ can be represented as an intersection of DVR's which are obtained by localizing at  its height-one primes,  i.e.,  $ B = \bigcap \{ B_P \mid P \in \sss(B), \hgt(P)$ $= 1\}.$
 It is well known that this collection is locally finite.  Now, note that $A$ is a local domain with maximal ideal   $N_A := T_3k[T_3]_{(T_3)} + (T_1, T_2) k(T_3)[T_1, T_2]_{(T_1, T_2)}$.  
 Choose any valuation overring $W$ of $A$ such that $M_W$ (the maximal ideal of $W$) is generated by
   $T_3$ and lies over the maximal ideal of $A$.  It is easy to see that many such valuation domains of the field $k(T_1, T_2, T_3)$ exist (e.g., let $W'$ be a valuation overring of $B$ with maximal ideal $M'$ lying over $M_B$ and such that the residue field $W'/M'$ is canonically isomorphic to   $k(T_3)$, which is the residue field $B/M_B$,  then the domain 
 $V + M'$, with $V$ as in the previous paragraph,  can serve as the desired domain $W$).   Now, the intersection
 $ R := \bigcap \{ B_P \mid P \in \sss(B), \hgt(P)= 1\} \bigcap W $
 is clearly a locally finite intersection.  We claim that any choice, as above, of the domain $W$ will yield $R = A$.   
 
 To prove our claim, we note first that it is obvious that $R$ is an overring of $A$.  So, we need to prove that $R \subseteq A$.  Observe that the ideal $M_B$ is an ideal of $A$ as well as of $B$.  It follows easily that $M_B$ is a prime ideal of $R$,  since $R \subset B$. Then, given an element $r \in R$, we can write
$r = \psi + f$, where $\psi\in  k(T_3)$ and $f \in M_B$.  However, $f \in M_B \subseteq W$  and so $\psi \in W$.  It is clear though that $W \cap  k(T_3) = V$.  It follows that $\psi \in V$ and so  $r \in A$.  Hence, we have proven that $R \subseteq A$.  
 \end{ex}

The following Proposition is the key step in proving the main results of the section.

\begin{prop}\label{prufer}
Let $A$ be a Pr\"ufer domain and $K$ be the quotient field of $A$. Let $Y$ be a subset of $Z:=\zar(K| A)$ such that $A=\bigcap\{ V \mid V\in Y\}$, and let $\gamma:\zar(K| A)\longrightarrow \sss(A)$, be the canonical map (defined by sending a valuation domain $V\in \zar(K|A)$ into its center in $A$). Then, 
$
\gamma^{-1}(\mmm(A))\subseteq \ad^{\mbox{\tiny{\texttt{cons}}}}(Y)
$.

\end{prop}
\noindent
\textsc{Proof.} Let  $M$ be  a maximal ideal of $A$. Since $A$ is a Pr\"ufer domain, the $t$-operation on $A$ coincides with the identity \cite[Theorem 22.1(3)]{Gilmer},  thus obviously $M$ is a $t-$maximal $t-$ideal of $A$.
 Now, we are able to apply \cite[Proposition 2.8(ii)]{calota} and, so,  there exists an ultrafilter $\ms U \in \beta(Y)$ such that 
$$M=   \{x\in A\mid \gamma^{-1}(V(x))\cap Y\in \ms U\}.$$
On the other hand,  the collection of sets
$$
\ms{V} := \{X' \subseteq \gamma(Y) \mid \gamma^{-1}(X') \cap Y\in \ms U\}
$$ 
is an ultrafilter  on $\gamma(Y)$ (precisely, with the notation of \cite[Lemma 2.1(4)]{fi-fo-lo0},   
 $\ms{V} ={ \ms{U}}_{\!{\gamma}}$, where for simplicity we have still denoted by $\gamma$ the map $\gamma|_Y : Y \rightarrow \gamma(Y)$) and, moreover,
$$
\begin{array}{rl}
P_{\ms V} :=& \{ x \in A \mid V(x) \cap \gamma(Y) \in \ms V \} \\
=&\{ x \in A \mid  \gamma^{-1}(V(x) \cap \gamma(Y) ) \cap Y \in \ms U  \} \\
=&\{ x \in A \mid  \gamma^{-1}(V(x)) \cap Y \in \ms U \} = M\,.
\end{array}
$$
 Moreover, if $A_{\ms U}$ is the ultrafilter limit valuation domain of $K$ associa\-ted to ${\ms U}\in \beta(Y)$ (i.e., $A_{\ms U} = \{x\in K \mid B_x\cap Y\in \ms U\}$, Corollary \ref{ultrafilter-zar}), then by \cite[Proposition 2.10(i)]{calota},  we have  $\gamma(A_{\ms U})= P_{\ms V} = M$. Therefore,  $\mmm(A)\subseteq \ad^{\mbox{\tiny{\texttt{cons}}}}(\gamma(Y))$.  Since $\gamma$ is continuous and closed   with respect to the  constructible topology   \cite[Theorem 3.8]{fi-fo-lo0}, it follows that 
$\ad^{\mbox{\tiny{\texttt{cons}}}}(\gamma(Y))=\gamma(\ad^{\mbox{\tiny{\texttt{cons}}}}(Y))$. 

 Moreover,  since $A$ is a Pr\"ufer domain, by \cite[Proposition 2.2]{dofo} $\gamma$ is also injective and, hence, $\gamma^{-1}(\mmm(A))\subseteq \ad^{\mbox{\tiny{\texttt{cons}}}}(Y)$. \hfill$\Box$
\medskip

Let $A$ be a domain, $K$ be the quotient field of $A$, and $T$  be an indeterminate over $K$. For each subset $Y$ of $Z:=\zar(K|A)$, we   set 
$$
Y_0:=\{V(T) \mid V\in Y\}, \qquad  Y^{\uparrow}:=\{V\in Z:V\supseteq W,\mbox{ for some }W\in Y\},
$$ 
where  $Y^\uparrow$ coincides with $Y^{\mbox{\tiny{\texttt{gen}}}}$ \emph{the Zariski--generic closure of $Y$}, i.e., the generic closure of $Y$ in $Z$, with respect to the Zariski topology}  (since $V\preccurlyeq  W  \mbox{ (in $Z^{\mbox{\tiny{\texttt{zar}}}}$)}\, :\Leftrightarrow \, V\supseteq W$).
{ Recall that, in Corollary \ref{cisiamoquasi2}(1), we introduced a general form of the Kronecker function ring, by setting $
\kr(K|A) = \bigcap \{ V(T) \mid V \in Z \} = \bigcap Z_0 =: \kr(Z)$. Now, we can extend this notion for $Y \subseteq Z$, by setting
$$
\kr(Y):= \bigcap Y_0 = \bigcap \{ V(T) \mid V \in Y \}\,
$$
 which is called the {\it $K$--function ring associated to $Y$}}.
  We recall that an integrally closed domain $A$ is a {{\it vacant domain} if, for each  $Y \subseteq Z$ such that $A = \bigcap Y$, then $\kr(Y) =\kr(Z)$  \cite[Definition 2.1.11]{Fabbri}.}

\begin{thm}\label{prufer2}
Let $K$ be a field and $C$ a closed subset of $\zar(K)^{\mbox{\tiny{\texttt{cons}}}}$.    Let $ (C^{\uparrow})_0 =\{W(T)\mid W\in C^{\uparrow}\}$.  Then, $
\zar(K(T)|\kr(C))=(C^{\uparrow})_0$.
\end{thm}
\noindent
\textsc{Proof.}
 The inclusion $\supseteq $ is obvious.  For the converse, let $\widetilde{W} \in \zar(K(T)|\kr(C))$. By Proposition \ref{cruciale}, we can suppose that  $ \widetilde{W} = W(T)$, for some $W\in \zar(K)$.
We want to show that $W\supseteq V$, for some $V\in C$.  Let $\varphi: \zar(K(T)) \rightarrow \zar(K)$ be the canonical map (Corollary \ref{embedding}). 
Since $\varphi|_{\zar_0(K(T))}: \zar_0(K(T))^{\mbox{\tiny{\texttt{cons}}}} \rightarrow \zar(K)^{\mbox{\tiny{\texttt{cons}}}}$ is a homeomorphism (Proposition \ref{ultrahomeo}),
  then the set 
$$
({\varphi|_{\zar_0(K(T))}})^{-1}(C)=\{V(T)\mid V\in C\}=C_0
$$ 
is closed both in $\zar_0(K(T))^{\mbox{\tiny{\texttt{cons}}}}$ and $\zar(K(T))^{\mbox{\tiny{\texttt{cons}}}}$ (Remark  \ref{infinity}).  Consider the natural map  $\gamma\!: \zar(K(T)|\kr(C))^{\mbox{\tiny{\texttt{cons}}}} \rightarrow \sss(\kr(C))^{\mbox{\tiny{\texttt{cons}}}}$\!, defined by sending a valuation overring of $\kr(C)$ into its center on $\kr(C)$.   Since the Kronecker function ring  $\kr(C)$ is, in particular, a Pr\"ufer domain with quotient field $K(T)$ (Proposition \ref{basic}(3))  then, from Proposition \ref{prufer},  it follows immedia\-tely that $\gamma^{-1}(\mmm(\kr(C)))\subseteq C_0$. 
Set $A(C):= \bigcap\{V \mid V\in C\}$. Now, by Zorn's Lemma, we can find a minimal valuation overring of $\kr(C)$ which, by Proposition \ref{cruciale},   is of the form $V'(T)$, for some $V'\in \zar(K|A(C))$, such that $W(T) \supseteq V'(T)$.  Then, by applying \cite[Corollary 19.7]{Gilmer} { (and, again, Proposition \ref{cruciale}),} we have  $\zar_{\rm min}(\kr(C))\subseteq \gamma^{-1}(\mmm(\kr(C)))$. Since,  by what we observed above, $ \gamma^{-1}(\mmm(\kr(C))) \subseteq C_0$, then   $V'(T)\in C_0$.\hfill$\Box$

\begin{oss} \label{4:6} \rm  Note that, with the notation and assumptions of Theo\-rem \ref{prufer2}, 
$$(C^{\uparrow})_0 = (C_0)^{\uparrow} := \{\widetilde{W}\in \zar(K(T)) \mid  \widetilde{W} \supseteq V(T), \mbox{ for some } V \in C\}\,.$$
As a matter of fact, $\{\widetilde{W}\in \zar(K(T)) \mid  \widetilde{W} \supseteq V(T), \mbox{ for some } V \in C\}  = $ $\{\widetilde{W}\in \zar(K(T)|\kr(C)) \mid  \widetilde{W} \supseteq V(T), \mbox{ for some } V \in C\}$.  By Proposition \ref{cruciale}, we have  $\zar(K(T)|\kr(C)) = \zar_0(K(T)|\kr(C))$, thus $(C_0)^{\uparrow} =  \{W(T) \in \zar(K(T)) \mid  W \in \zar(K),\, W(T) \supseteq V(T), \mbox{ for some } V \in C\} = (C^{\uparrow})_0$.
\end{oss}

\medskip

We  recall  some properties of semistar operations.  Let $A$ be an integral domain and let $K$ be the quotient field of $A$. We denote by $\FF(A)$ the set of all the nonzero $A$-submodules of $K$, and by $\f(A)$ the set of all the nonzero finitely generated $A$-submodules of $K$. 
 A map
$
\star:\FF(A)\rightarrow \FF(A), \, E\mapsto E^\star
$,
is called {\it a semistar operation on} $A$ if, for each $0\neq x\in K$ and for all $E, F\in \FF(A)$, the following properties hold:
($\boldsymbol{\star_1}$)\  $(xE)^\star=xE^\star$; \ 
($\boldsymbol{\star_2}$)\ $E\subseteq F\Rightarrow E^\star \subseteq F^\star$; \
($\boldsymbol{\star_3}$)\  $E\subseteq E^\star$ and $(E^\star)^\star=E^\star$.

A {\it semistar operation of finite type} $\star$ on $A$ is a semistar operation such that, for every $E\in\FF(A)$, 
$$
E^\star= E^{\stf} := \bigcup\{F^\star\mid F\in \f(A), \  F\subseteq E\}\,.$$

An {\it \texttt{e.a.b.} semistar operation} $\star$ on $A$ is a semistar operation such that, for all $F, G, H\in \f(A)$, $(FG)^\star\subseteq (FH)^\star$ implies $G^\star\subseteq H^\star$.

A valuation domain $V\in \zar(K|A)$ is called a {\it $\star$--valuation overring of} $A$ if $F^\star\subseteq FV$, for each $F\in \f(A)$. We denote by $\zar^\star(K|A)$ the collection of all the $\star$--valutation overring of $A$.

Important classes of examples of semistar operations are obtained as follows. Let  $\mathcal S$ be a nonempty family of overrings of $A$.  
Then the map  $\wedge_{\mathcal S}:\FF(A) \rightarrow \FF(A)$, 
$E\mapsto \bigcap\{ES\mid S\in\mathcal S\}$, 
defines a semistar operation on $A$ \cite[Theorem 1.2(C)]{fohu}. 
In particular, given a nonempty subset $Y$ of $\zar(K|A)$, the semistar operation $\wedge_Y$ is  \texttt{e.a.b.}, by \cite[Proposition 7]{folo2}. We say that a {\it semistar operation} $\star$ on $A$ is {\it complete} if $\star= b(\star):=\wedge_{\zar^\star(K|A)}$. For any semistar operation $\star$ on  $A$, it is easily seen that $F^{b(\star)}V=F^\star V$, for each $F\in \f(A)$ and $V\in \zar^\star(K|A)$. Thus,  $b(b(\star))=b(\star)$ and $b(\star)$ is a complete semistar operation.   The {\it $b$--operation} on $A$ is the \texttt{e.a.b.} semistar  operation  defined by
$b:= \wedge_{\zar(K|A)}$ and, obviously, $b \leq b(\star)$ (i.e., $E^b \subseteq E^{b(\star)}$ for each $E \in \FF(A)$) for all semistar operations $\star$ on $A$.

\medskip

\begin{oss}\label{generic}
Let $K$ be a field and $A$ be a subring of $K$. If $Y$ is a nonempty subset of $\zar(K|A)$, then it is immediately seen that $\wedge_Y=\wedge_{Y^\uparrow}$.  Moreover, if $Y$ is a  quasi-compact subset of  $\zar(K|A)^{\mbox{\tiny{\texttt{zar}}}}$, the subset $Y_{\mbox{\tiny{\texttt{min}}}}$, consisting  of the minimal elements of $Y$, is nonempty and it is easy to see that $\wedge_{Y^\uparrow}= \wedge_Y= \wedge_{Y_{\mbox{\Tiny{\texttt{min}}}}}$.
\end{oss}

\medskip

Let $T$ be an indeterminate over $A$ and $f\in A[T]$. We shall denote by $\cc(f)$ the content of the polynomial $f$. If $\star$ is an \texttt{e.a.b.} semistar operation on $A$, $$\kr(A, \star) := \{f/g \in K(T) \mid f,g\in A[T],\  g\neq 0,  \mbox{ and }\ \cc(f)^\star\subseteq \cc(g)^\star \}
$$
is called the {\it $\star$--Kronecker function ring of $A$}.
It is well known that $\kr(A, \star)$ 
is a B\'ezout domain with quotient field $K(T)$ \cite[Theorems 5.1 and 3.11(3)]{folo1}.
Note that, if   $Y$ is a subset of $\zar(K|A)$, by \cite[Corollary 3.8]{folo1}, the $\wedge_Y$--Kronecker function ring of $A$, $\kr(A, \wedge_Y)$, coincides with the $K$--function ring  $\kr(Y)$ (introduced just before Theorem \ref{prufer2}). 

\medskip

Now, we give an application of the ultrafilter topology for characte\-ri\-zing when two \texttt{e.a.b.} semistar operations of finite type are equal.

\begin{thm}\label{semistarclo}
Let $A$ be an integral domain with quotient field $K$ and $Y', Y''\subseteq \zar(K|A)$. Then, the following conditions are equivalent.
\begin{enumerate}
\item[\rm (i)] $(\wedge_{Y'})_f=(\wedge_{Y''})_f$.
\item[\rm (ii)] The sets $\ad^{\mbox{\tiny{\texttt{cons}}}}(Y'),\ad^{\mbox{\tiny{\texttt{cons}}}}(Y'')$ have the same Zariski--generic closure in $\zar(K|A)$, i.e., $\ad^{\mbox{\tiny{\texttt{cons}}}}(Y')^\uparrow=\ad^{\mbox{\tiny{\texttt{cons}}}}(Y'')^\uparrow$.
\end{enumerate}
\end{thm}
\noindent
\textsc{Proof.} Let $T$ be an indeterminate over $K$. By \cite[Remark 3.5(b)]{folo1}, it is enough to show that condition (ii) is equivalent to the following   
\begin{enumerate}
\item[\rm (i$'$)]   $\kr(A, {\wedge_{Y'}})=\kr(A, {\wedge_{Y''}})$.
\end{enumerate}

(ii)$\Rightarrow$(i$'$). Assume that the equality  $\ad^{\mbox{\tiny{\texttt{cons}}}}(Y')^\uparrow=\ad^{\mbox{\tiny
{\texttt{cons}}}}(Y'')^\uparrow$ holds. Keeping in mind the notation introduced before Theorem \ref{prufer2} and applying Corollary \ref{cisiamoquasi3}(1), it follows easily that, inside $\zar(K(T))$, $\ad^{\mbox{\tiny{\texttt{cons}}}}(Y'_0)^\uparrow=\ad^{\mbox{\tiny
{\texttt{cons}}}}(Y''_0)^\uparrow$.   By using Proposition \ref{closure-intersection} and Remark \ref{generic}, we have
$$
\begin{array}{rl}
\bigcap Y'_0= &\hskip -4pt \bigcap \ad^{\mbox{\tiny{\texttt{cons}}}}(Y'_0)= \bigcap \ad^{\mbox{\tiny{\texttt{cons}}}}(Y'_0)^\uparrow =\bigcap \ad^{\mbox{\tiny{\texttt{cons}}}}(Y''_0)^\uparrow \\
=&\hskip -4pt \bigcap \ad^{\mbox{\tiny{\texttt{cons}}}}(Y''_0)=\bigcap Y''_0,
\end{array}
$$
and thus $\kr(A, \wedge_{Y'})=\kr(A, \wedge_{Y''})$, in view of \cite[Corollary 3.8]{folo1}.

(i$'$)$\Rightarrow$(ii). Set $B:=\kr(A, \wedge_{Y'})=\kr(A, \wedge_{Y''})$. 
 By using \cite[Corollary 3.8]{folo1}, Proposition \ref{closure-intersection}, Theorem \ref{prufer2}  and Remark \ref{4:6}, it follows that 
$$\ad^{\mbox{\tiny{\texttt{cons}}}}(Y'_0)^\uparrow= \zar(K(T)|B)=
\ad^{\mbox{\tiny{\texttt{cons}}}}(Y''_0)^\uparrow,$$
and thus the conclusion is clear, again by Corollary \ref{cisiamoquasi3}(1). \hfill$\Box$

\begin{cor}\label{vacantcharacterization}
Let $A$ be an integrally closed domain. Then, the following conditions are equivalent.
\begin{enumerate}
\item[\rm(i)] $A$ is a vacant domain.
\item[\rm(ii)] For each representation $Y\subseteq \zar(K|A)$ of $A$ (i.e., $\bigcap Y=A$), we have $\ad^{\mbox{\tiny{\texttt{cons}}}}(Y)^\uparrow=\zar(K|A)$.
\end{enumerate}
\end{cor}
\noindent
\textsc{Proof.} { Set $Z:= \zar(K|A)$.} 

(i)$\Rightarrow$(ii). Assume $A$ vacant and take a subset $Y\subseteq Z$ such that $\bigcap Y=A$.  By \cite[Proposition 3.3]{folo1}, we have $\kr(A, \wedge_Y)=\kr(Y)={ \kr(Z)=
\kr(A, b)}$, and thus 
$$
(\wedge_Y)_f=b={ \wedge_{Z}
= {(\wedge_{Z})}_f.}
$$
{ The conclusion follows immediately from} Theorem \ref{semistarclo}.

(ii)$\Rightarrow$(i). Take a subset $Y$ of $\zar(K|A)$ such that $\bigcap Y=A$. By assumption and Theorem \ref{semistarclo}, it follows that $(\wedge_Y)_f={ (\wedge_Z)_f = \wedge_Z} =  b$, and thus $\kr(Y) = {  \kr(A, \wedge_Y)= \kr(A, \wedge_Z) }=\kr(Z)$. This proves that $A$ is vacant. \hfill$\Box$

\medskip

From the previous theorem, we deduce immediately the following 
\begin{cor}\label{densevacant}
Let $A$ be an integrally closed domain. If each representation of $A$ is dense in $\zar(K|A)^{\mbox{\tiny{\texttt{cons}}}}$, then $A$ is a vacant domain. 
\end{cor}

\begin{ex}
Let $K$ be a field and  $T_1, T_2$ two indeterminates over $K$. Consider the pseudo-valuation domain  
 $A:= K +T_2K(T_1)[T_2]_{(T_2)}$ 
with associated valuation domain   $V:=K(T_1)[T_2]_{(T_2)}$ of  $K(T_1, T_2)$. Let $p:V\rightarrow  K(T_1)$ be the canonical projection of $V$ onto its residue field  $K(T_1)$ and so $A=p^{-1}(K)$.    Then, by \cite[Exercise 12, page 409]{Gilmer}, the domain $A$  is a vacant domain. It is easily seen that the set $C:=\{p^{-1}(W')\mid W'\in  \zar(K(T_1)|K)\} \subset \zar(K(T_1, T_2)|A)$ is a representation of $A$, and that it is closed,  with respect to the constructible topology of  $\zar(K(T_1, T_2|A)$, since $C=\{W\in \zar(K(T_1, T_2)|A)\mid W\subseteq V\}    =\bigcap_{z \in  K(T_1, T_2) \setminus V} (\zar(K(T_1,T_2)|A) \setminus B_z)=\ad^{\mbox{\tiny \texttt{zar}}}(\{V\})$. Thus, the converse of the previous Corollary \ref{densevacant}  does not hold in general.

Note that this example shows also that, in the statement of Theorem \ref{prufer2}, we need to consider $C^\uparrow$ and not just $C$, since  in this case $ \zar(K(T; T_1, T_2)|\kr(A, \wedge_C)) =\zar(K(T; T_1, T_2)|\kr(C))=(C^{\uparrow})_0 \supsetneq C_0$.
\end{ex}


Now, we prove that the property of being ``complete'' for a semistar operation can be caracterized by a ``compactness'' property for a suitable subspace of the Zariski-Riemann surface.
 
\begin{thm}\label{complete-compact}
Let $A$ be an integral domain with quotient field $K$ and $\star$ be a semistar operation on $A$. Then, the following conditions are equivalent.
\begin{enumerate}
\item[\rm (i)] $\star$ is \texttt{e.a.b.} of finite type.

\item[\rm (ii)] $\star$ is complete.

\item[\rm { (iii)}] There exists a closed subset { $Y$} of $\zar(K|A)^{\mbox{\tiny{\texttt{cons}}}}$  such that {$Y = {Y}^\uparrow$ and $\star=\wedge_{Y}$.}

\item[\rm { (iv)}] There exists a compact  subspace { $Y'$} in $\zar(K|A)^{\mbox{\tiny{\texttt{cons}}}}$  such that { $\star=\wedge_{Y'}$}.

\item[\rm (v)] There exists a { quasi-}compact { subspace of} $Y''$ of $\zar(K|A)^{\mbox{\tiny{\texttt{zar}}}}$  such that  $\star=\wedge_{Y''}$.
\end{enumerate}
\end{thm}
\noindent
\textsc{Proof.} (i)$\Leftrightarrow$(ii) depends on the fact that if $\star$ is \texttt{e.a.b.}, then $\stf = b(\star)$ \cite[Proposition 9]{folo2}.

 Let $T$ be an indeterminate on $K$ and let $\varphi:\zar(K(T))\longrightarrow \zar(K)$  be the canonical surjective map, defined in Proposition \ref{embedding}. 

(ii)$\Rightarrow$(iv). By \cite[Theorem 3.5]{folo4}, $\zar^{\star}(K|A)=\varphi(\zar(K(T)|
\kr(A, \star)))$, and thus (by Proposition \ref{ultrahomeo}) it is closed in $\zar(K|A)^{\mbox{\tiny{\texttt{cons}}}}$ or, equivalently, compact, in the compact Hausdorff space $\zar(K|A)^{\mbox{\tiny{\texttt{cons}}}}$.  Then,  the conclusion follows by taking ${ Y'}:=\zar^\star(K|A)$ (since, by definition,  $b(\star)=\wedge_{\zar^\star(K|A)}$). 

(iv)$\Rightarrow$(ii).  As observed above the compact subspaces of $\zar(K|A)^{\mbox{\tiny{\texttt{cons}}}}$ are exactly the closed subsets. Take a closed set { $Y'$} of  $\zar(K|A)^{\mbox{\tiny{\texttt{cons}}}}$ such that  $\star = \wedge_{Y'}$.  Set 
$
{Y'}^\uparrow:=\{W\in \zar(K|A)\mid W\supseteq V, \mbox{ for some } V\in Y' \}$.
 By \cite[Corollary 3.8]{folo1}, we have  
$
\kr(A, \wedge_{ Y'})= \bigcap\{V(T)\mid V\in { Y'}\}=:\kr({ Y'})
$.
On the other hand, since $Y'$ is closed, by Theorem \ref{prufer2}, it follows that $\zar(K(T)|\kr(A, {\wedge_{ Y'}}))=({ Y'}^\uparrow)_0 =
\{W(T)\mid W\in { Y'}^\uparrow\}$.
Therefore, as above (by \cite[Theorem 3.5]{folo4}), $\zar^{\wedge_{Y}}(K|A) = \varphi(\zar(K(T)|\kr(A, {\wedge_Y}))) = Y^\uparrow$ and so $Y^\uparrow$ is also a closed subspace of $\zar(K|A)^{\mbox{\tiny{\texttt{cons}}}}$.   Since by definition $b({\wedge_Y})= \wedge_{\zar^{\wedge_{Y}}(K|A)}= \wedge_{Y^\uparrow}$, then the conclusion is immediate, by Remark \ref{generic}. 

(iii)$\Rightarrow$(iv) is trivial since, as observed above, closed coincides with compact in $\zar(K|A)^{\mbox{\tiny{\texttt{cons}}}}$.

(iv)$\Rightarrow$(v) is obvious, by  \cite[Theorem 3.4(1)]{fi-fo-lo0}). 

(v)$\Rightarrow${ (iii)}. Take a set $Y''$ as stated in (v). Then,  (iii)  follows immediately from  Proposition \ref{genericultra-spectral} and Remark \ref{generic}, by taking $Y :=Y''^{\uparrow}$.
\hfill$\Box$

\smallskip

\begin{cor}\label{wedge-hat}
Let $A$ be an integral domain and $K$   its quotient field. Let $Y$ be a subset of $\zar(K|A)$ and set $\widehat Y:=\ad^{\mbox{\tiny{\texttt{cons}}}}(Y)^\uparrow$. Then,  $(\wedge_Y)_f=\wedge_{\widehat Y}
{
=\wedge_{
{\boldsymbol{\mbox{\it \tiny\texttt{Cl}}}}^
{\mbox{\footnotesize{\texttt{\Tiny{cons}}}}}(Y)
}
}$.
\end{cor}
\noindent
\textsc{Proof.} In view of Proposition \ref{genericultra-spectral}, $\widehat Y$ is closed, with respect to the constructible topology. Thus $\wedge_{\widehat Y}$ is of finite type, by Theorem \ref{complete-compact}, and hence the equali\-ty $(\wedge_Y)_f=\wedge_{\widehat Y}$ follows immediately by Theorem \ref{semistarclo}, since  $\ad^{\mbox{\tiny{\texttt{cons}}}}(Y)^\uparrow = {\widehat Y}^\uparrow (=\ad^{\mbox{\tiny{\texttt{cons}}}}(\widehat Y^\uparrow))$. { Moreover, the semistar operation $\wedge_{\boldsymbol{\mbox{\it \tiny\texttt{Cl}}}^
{\mathtt{cons}}(Y)}$ is of finite type, by Theorem \ref{complete-compact}, and thus the last equality follows by applying Theorem \ref{semistarclo}}. \hfill$\Box$

\bigskip

The next example illustrates the possibility that the sets $Y,\ Y'$ and $Y''$  in Theo\-rem \ref{complete-compact}  can form a proper chain of sets.

\begin{ex}  Let $k$ be a field and  $T_1, T_2$
 two indeterminates over $k$. Let $A$ be the two-dimensional,  integrally closed, local domain $k[T_1, T_2]_{(T_1, T_2)}$ 
 with quotient field $K:=  k(T_1, T_2)$.  
 Let $\star$ be the $b$--operation on $A$.  It is well known that the $b$--operation is an \texttt{e.a.b.} operation of finite type.  Hence, it satisfies the equivalent conditions of Theorem \ref{complete-compact}.  Our goal is to show that there is a great deal of flexibility in the choice of the sets $Y, Y'$ and $Y''$ in the theorem.  First, note that if the valuation domains in  $\zar(K|A) $  are ordered by inclusion then any chain is finite \cite[Corollary 30.10]{Gilmer} and, hence, { obviously} there are minimal elements.  Any such minimal valuation overring $V$ will necessarily have maximal ideal $M_V$ lying over  the maximal ideal  $(T_1, T_2)$ of $A$.  The standard definition of the $b$--operation involves extending an ideal  (or, more generally a sub-$A$-module of $K$)  to all valuation overrings.   It is clearly sufficient to extend to just those valuation overrings that are minimal.  So, any subcollection of  $\zar(K|A)$ which contains all the minimal elements will generated the $b$--operation.  It is not clear that the collection of minimal valuation overrings is closed under the Zariski or the constructible topology.

\begin{itemize}

\item  Consider the  members of  $\zar(K|A)$ which do not contain the elements  $\frac{1}{T_1} , \frac{1}{T_2}$.  This is a closed, quasi-compact subset of  $\zar(K|A)^{\mbox{\tiny{\texttt{zar}}}}$.  It can also be thought of as being those valuation domains in  $\zar(K|A)$ whose maximal ideal dominates  $(T_1, T_2)$ in $A$.   Hence, it contains the minimal valuation overrings and is sufficient to generate the $b$ operation.   We can let this collection be denoted by 
$Y''$ in Theorem 4.14.

\item The set $Y''$, described above, is { a (proper) closed subset of}   $\zar(K|A)^{\mbox{\tiny{\texttt{zar}}}}$.  Hence, it is also closed in  $\zar(K|A)^{\mbox{\tiny{\texttt{cons}}}}$.  Moreover, any closed subset of 
 $\zar(K|A)^{\mbox{\tiny{\texttt{cons}}}}$ is compact.  Hence, to obtain our set  $Y'$, we can choose any closed subset of  $\zar(K|A)^{\mbox{\tiny{\texttt{cons}}}}$ which contains 
 $Y''$.  Since any single point is closed in  $\zar(K|A)^{\mbox{\tiny{\texttt{cons}}}}$, we can let  $Y'$ be the union of $Y''$ and any other single valuation overring, for example, the localization of $A$ at a height-one prime.
 
 \item The set $Y$ should contain all overrings of its members.  An obvious choice then is to let $Y$ be all of  $\zar(K|A)^{\mbox{\tiny{\texttt{cons}}}}$.  Since this is the entire space it is trivially closed {(in $\zar(K|A)^{\mbox{\tiny{\texttt{cons}}}}$)} and generates the $b$--operation.

 \end{itemize}
 
 This then gives three different sets  $Y'' \subset Y' \subset Y$ 
 with the notation of Theorem \ref{complete-compact}, all associated with the same (semi)star operation.

\end{ex}

By using Remark \ref{inv-spectral}, we can restate Corollaries  \ref{vacantcharacterization} and \ref{wedge-hat} as follows:

\begin{cor}
Let $A$ be an integrally closed domains and $K$ be its quotient field. Then the following conditions are equivalent.
\begin{enumerate}
\item[\rm (i)] $A$ is a vacant domain.
\item[\rm (ii)] Each representation of $A$ is dense in $\zar(K|A)$, with respect to the inverse
 topology.
\end{enumerate}
\end{cor}

\begin{cor}
Let $A$ be an integral domain and $K$   its quotient field. Let $Y$ be a subset of $\zar(K|A)$. Then,  $(\wedge_Y)_f=\wedge_{{\tiny{ \ad}}^{\mbox{\Tiny{\texttt{\emph{inv}}}}}(Y)}$.
\end{cor}
\bigskip
\noindent {\bf Acknowledgement}

\noindent The authors would like to thank the referee for suggestions which substantially improved the paper.

 

\end{document}